\documentclass[11pt,french]{article}
\usepackage{latexsym}
\usepackage[all]{xy}
\usepackage{amsmath}
\usepackage{amssymb}
\usepackage{amsfonts}

\newtheorem{thm}{Th\'eor\`eme}[section]
\newtheorem{prop}[thm]{Proposition}

\newtheorem{df}[thm]{D\'efinition}
\newtheorem{cor}[thm]{Corollaire}

\begin{document}

\title{\textbf{Alg\`ebres simpliciales $S^{1}$-\'equivariantes et th\'eorie de de Rham}}
\bigskip
\bigskip

\author{\bigskip\\
\textbf{Bertrand To\"en} \\
\small{Laboratoire Emile Picard} \\
\small{Universit\'{e} Paul Sabatier, Bat 1R2}
\small{31062 Toulouse Cedex 9}\\
\small{France}\\
\medskip \\
\textbf{Gabriele Vezzosi}\\
\small{Dipartimento di Matematica Applicata
``G. Sansone''} \\ 
\small{Universit\`a di Firenze }\\
\small{Italy}}

\date{Avril 2009}

\maketitle

\begin{abstract}
Ce travail accompagne \cite{tv} et a pour objectif de fournir de plus
amples d\'etails sur la comparaison entre fonctions sur les espaces
des lacets d\'eriv\'es et th\'eorie de de Rham. Pour une
$k$-alg\`ebre commutative $A$, lisse sur $k$ de caract\'eristique nulle, 
nous montrons que deux objets $S^{1}\otimes A$ et $\epsilon(A)$
se d\'eterminent mutuellement, et ce fonctoriellement en $A$.
L'objet $S^{1}\otimes A$ est la $k$-alg\`ebre simpliciale
$S^{1}$-\'equivariante obtenue en tensorisant $A$ par 
le groupe simplicial $S^{1}:=B\mathbb{Z}$. L'objet 
$\epsilon(A)$ est l'alg\`ebre de de Rham de $A$, munie
de la diff\'erentielle de de Rham et consid\'er\'ee
comme une \emph{$\epsilon$-dg-alg\`ebre} (i.e. une
alg\`ebre dans une certaine cat\'egorie mono\"\i dale de $k[\epsilon]$-dg-modules, 
o\`u $k[\epsilon]:=H_{*}(S^{1},k)$). Nous construisons une
\'equivalence $\phi$, entre la th\'eorie homotopique des $k$-alg\`ebres simpliciales
$S^{1}$-\'equivariantes et celle des
$\epsilon$-dg-alg\`ebres, et nous montrons l'existence d'une
\'equivalence fonctorielle $\phi(S^{1}\otimes A)\sim \epsilon(A)$.
Nous d\'eduisons de cela la comparaison annonc\'ee,  identifiant les fonctions 
$S^{1}$-\'equivariantes sur $LX$, l'espace des lacets d\'eriv\'e
d'un $k$-sch\'ema $X$ lisse, et la cohomologie de de Rham alg\'ebrique
de $X/k$. Nous d\'eduisons aussi des versions fonctorielles et multiplicatives
des th\'eor\`emes de d\'ecomposition de l'homologie de Hochschild (du type HKR), pour des
$k$-sch\'emas s\'epar\'es quelconques.
\end{abstract}

\section*{Introduction}

Dans ce travail nous comparons la th\'eorie des fonctions sur l'espace
des lacets d\'eriv\'es d'un sch\'ema de caract\'eristique nulle $X$, au sens de la g\'eom\'etrie
alg\'ebrique d\'eriv\'ee (voir \cite{to,tv',tv}), avec sa th\'eorie de 
de Rham. Cette comparaison est annonc\'ee dans \cite{bena} ainsi que dans 
notre r\'ecent travail \cite{tv}, et semble d\'ecouler d'une comparaison plus g\'en\'erale, mais
encore conjecturale, 
entre fonctions sur les espaces de lacets d\'eriv\'es et homologie cyclique. Dans ce travail
nous \'etablirons cette comparaison avec la th\'eorie de de Rham de mani\`ere directe, sans
avoir \`a passer par l'homologie cyclique. Une 
cons\'equence int\'eressante de cette approche directe est de fournir de nouvelles preuves, et une
nouvelle compr\'ehension, des th\'eor\`emes HKR pour les sch\'emas (dans le style de \cite{ye,sch}).

Soit $k$ un anneau commutatif de caract\'eristique nulle et $A$ une $k$-alg\`ebre
lisse sur $k$. Le sch\'ema d\'eriv\'e des lacets de $X=\mathrm{Spec}\, A$ est par
d\'efinition $LX:=\mathbb{R}Map(S^{1},X)$, o\`u $S^{1}:=B\mathbb{Z}$ est consid\'er\'e
comme un groupe simplicial (voir \cite{bena,to,tv} pour les d\'etails). Au niveau 
des anneaux de fonctions $LX$ est le spectre de la $k$-alg\`ebre
commutative simpliciale $S^{1}\otimes A$, o\`u l'on utilise ici l'enrichissement
simplicial naturel de la cat\'egorie des $k$-alg\`ebres simpliciales commutatives. Ainsi, 
le probl\`eme qui consiste \`a comparer les fonctions sur $LX$ et la th\'eorie
de de Rham de $X$ se r\'esume de mani\`ere purement alg\'ebrique \`a \'etudier
les relations entre $S^{1}\otimes A$ et l'alg\`ebre de de Rham de $A/k$. 
Le r\'esultat principal de ce travail affirme que la donn\'ee
de $S^{1}\otimes A$, muni de son action naturelle de $S^{1}$, est \'equivalente, \`a
homotopie pr\`es, 
\`a la donn\'ee de l'alg\`ebre de de Rham de $A/k$, munie de sa diff\'erentielle de
de Rham.

Afin de donner un \'enonc\'e plus pr\'ecis 
consid\'erons $DR(A)=Sym_{A}(\Omega_{A/k}^{1}[1])$ l'alg\`ebre
de de Rham de $A/k$, pour laquelle $\Omega_{A}^{n}$ est plac\'e en degr\'e
$-n$. Muni  de la diff\'erentielle nulle $DR(A)$ devient une
$k$-dg-alg\`ebre commutative (\emph{cdga} pour faire court). De plus, la diff\'erentielle de
de Rham muni cette cdga d'une structure additionelle, \`a savoir celle
d'une $\epsilon-cdga$, c'est \`a dire d'une op\'eration de degr\'e -1, $\epsilon : DR(A) \longrightarrow
DR(A)[1]$, satisfaisant la r\`egle de Liebniz (au sens gradu\'e) par rapport
\`a la multiplication dans $DR(A)$ (voir \S 1 pour la d\'efinition pr\'ecise 
de $\epsilon-cdga$). Cette $\epsilon-cdga$ sera not\'ee $\epsilon(A)$.
D'autre part, on peut aussi former $S^{1}\otimes A$, la $k$-alg\`ebre
commutative simpliciale obtenue en tensorisant $A$ par le groupe simplicial $S^{1}=B\mathbb{Z}$.
L'objet $S^{1}\otimes A$ est naturellement muni d'une action du groupe simplicial 
$S^{1}$, qui op\`ere sur lui-m\^eme par translations, et est donc un objet
de $S^{1}-sk-CAlg$, la cat\'egorie des $k$-alg\`ebres simpliciales commutatives
$S^{1}$-\'equivariantes. Notre th\'eor\`eme principal (voir \ref{t1}) peut alors 
s'\'enoncer de la fa\c{c}on suivante. 

\begin{thm}\label{ti}
Il existe une \'equivalence, $\phi$, de la th\'eorie homotopique 
$S^{1}-sk-CAlg$ vers la th\'eorie homotopique $\epsilon-cdga$, qui est 
telle que pour tout $k$-alg\`ebre commutative lisse $A$ il existe un
isomorphisme fonctoriel
$$\epsilon(A) \simeq \phi(S^{1}\otimes A).$$
\end{thm}

La notion d'\emph{\'equivalence de th\'eorie homotopique} utilis\'ee dans
le th\'eor\`eme pr\'ec\'edent sera pour nous une \'equivalence de d\'erivateurs
au sens de Grothendieck (voir \cite{gr}). Ainsi, le th\'eor\`eme \ref{ti}
nous dit que pour toute petite cat\'egorie $I$, il existe une 
\'equivalence de cat\'egories homotopiques de diagrammes
$$\phi_{I} : Ho(S^{1}-sk-CAlg^{I}) \longrightarrow Ho(\epsilon-cdga^{I}),$$
fonctorielle en $I$. De plus, pour tout $I$-diagramme de $k$-alg\`ebres
commutatives lisses $A$, il existe un isomorphisme dans $Ho(\epsilon-cdga^{I})$
$$\phi(S^{1}\otimes A) \simeq \epsilon(A),$$
qui est non seulement fonctoriel en $A$, mais aussi en $I$. Le th\'eor\`eme
\ref{ti} sera en r\'ealit\'e une cons\'equence d'un r\'esultat plus g\'en\'eral, 
valable pour toute $k$-alg\`ebre commutative simpliciale $A$, affirmant l'existence
d'une \'equivalence
$$\mathbb{L}\epsilon(N(A)) \simeq \phi(S^{1}\otimes^{\mathbb{L}} A),$$
o\`u $N(A)$ est la $cdga$ obtenue par normalisation \`a partir de $A$, et 
$\mathbb{L}\epsilon$ et $S^{1}\otimes^{\mathbb{L}}-$ sont des versions 
d\'eriv\'ees des constructions $\epsilon$ et $S^{1}\otimes -$. \\

Un corollaire important, et imm\'ediat, de notre th\'eor\`eme principal est la
version suivante du th\'eor\`eme HKR.

\begin{cor}\label{ci}
Soit $X$ un $k$-sch\'ema s\'epar\'e. Alors, il existe un isomorphisme naturel
dans la cat\'egorie homotopique des faisceaux de $\mathcal{O}_{X}$-dg-alg\`ebres
commutatives
$$Sym_{\mathcal{O}_{X}}(\mathbb{L}_{X/k}[1])\simeq \mathcal{O}_{X}\otimes^{\mathbb{L}}_{\mathcal{O}_{X}\otimes^{\mathbb{L}}_{k}
\mathcal{O}_{X}}\mathcal{O}_{X},$$
o\`u $\mathbb{L}_{X/k}$ est le complexe cotangent de \cite{il}.
En particulier, si $X$ est lisse sur $k$ on a
$$Sym_{\mathcal{O}_{X}}(\Omega^{1}_{X/k}[1])\simeq \mathcal{O}_{X}\otimes^{\mathbb{L}}_{\mathcal{O}_{X}\otimes^{\mathbb{L}}_{k}
\mathcal{O}_{X}}\mathcal{O}_{X}.$$
\end{cor}

Il faut noter ici que ce corollaire redonne les r\'esultats 
de \cite{ye} et \cite{sch}. Il les am\'eliore aussi car nos isomorphismes
sont \textit{multiplicatifs}. Il n'est d'ailleurs pas tout \`a fait clair que 
les isomorphismes de \ref{ci} soient les m\^emes que ceux
de \cite{sch,ye}. 

Pour finir cette introduction, deux mots concernant la strat\'egie de la preuve
du th\'eor\`eme \ref{ti} et les difficult\'es que nous avons rencontr\'ees. 
Il faut en fait remarquer que le point crucial
est la construction de l'\'equivalence $\phi$. En effet, par d\'efinition, 
$S^{1}\otimes A$ est \emph{la $k$-alg\`ebre simpliciale $S^{1}$-\'equivariante
libre sur $A$}. De m\^eme, nous montrons (voir proposition \ref{p2}) que
$\epsilon(A)$ est \emph{la $\epsilon-cdga$ libre sur $A$}. Ainsi, une fois
l'\'equivalence $\phi$ construite on d\'eduit l'existence d'un isomorphisme
naturel $\epsilon(A)\simeq \phi(S^{1}\otimes A)$ formellement, par propri\'et\'e
universelle de ces deux objets: il suffit en effet que $\phi$ soit compatible
avec certains foncteurs qui oubient d'une part l'action de $S^{1}$ et 
d'autre part la $\epsilon$-structure. La construction d'une telle 
\'equivalence $\phi$ est donn\'ee dans notre \S 2 et se r\'ev\`ele plus
compliqu\'ee que nous le croyions d'abord. En effet, il se trouve que nous
n'avons pas trouv\'e d'approches directes reliant les th\'eories homotopiques
$S^{1}-sk-CAlg$ et $\epsilon-cdga$, et notre construction de $\phi$ 
passe par une chaine relativement longue d'\'equivalences de Quillen entre
plusieurs cat\'egories de mod\`eles auxiliaires. C'est pour cette raison que nous avons
choisi de formuler cette construction dans le contexte
des d\'erivateurs de Grothendieck, qui est relativement efficace afin de ne pas avoir
\`a trainer d'interminables chaines d'\'equivalences fonctorielles. Il est possible cependant, qu'un
approche plus directe existe, et nous pr\'etendons simplement ne pas l'avoir vue. 
A ce sujet, nous faisons remarquer qu'il existe d'une part une \'equivalence
de Quillen $N : sk-CAlg \longrightarrow cdga$, entre $k$-alg\`ebres simpliciales commutatives
et $cdga$, induite par le foncteur de normalisation de la correspondance de Dold-Kan
(voir \cite{ss}). D'autre part, il existe aussi une \'equivalence de Quillen 
$N : S^{1}-sk-Mod \longrightarrow k[\epsilon]-dg-mod$, entre les 
$k$-modules simpliciaux $S^{1}$-\'equivariants et les $k[\epsilon]$-dg-modules
(ici $k[\epsilon]=H_{*}(S^{1},k)$), qui elle aussi est induite par normalisation. 
Cependant, ces deux \'equivalences de Quillen ne semblent pas se promouvoir en une 
\'equivalence de Quillen $N : S^{1}-sk-CAlg \longrightarrow \epsilon-cdga$, induite
par le simple foncteur de normalisation. L'obstruction \`a l'existence d'un tel 
foncteur de normalisation provient du fait que $N : S^{1}-sk-Mod \longrightarrow k[\epsilon]-dg-mod$
n'est pas compatible avec les structures mono\"\i dales de ces deux cat\'egories
(qui sont induites par les produits tensoriels sur $k$ et utilisent donc
des structures de type co-alg\`ebre pour \^etre d\'efinies), 
m\^eme pas en un sens \emph{lax} de \cite{ss}. Le mophisme \emph{shulffle}
$N(E)\otimes N(F) \longrightarrow N(E\otimes F)$ n'est simplement pas
un morphisme de $k[\epsilon]$-dg-modules. Ainsi, 
pour $A \in S^{1}-sk-CAlg$, $N(A)$ poss\`ede bien une structure de $cdga$ d'une part, et une
structure de $k[\epsilon]$-dg-module d'autre part, mais ces deux structures ne v\'erifient 
pas les conditions de compatibilit\'e pour faire de $N(A)$ une $\epsilon-cdga$. 
Cela dit, il est int\'eressant de remarquer que notre \'equivalence $\phi$ sera probablement consid\'er\'ee
par les sp\'ecialistes comme \emph{bien connue}, et tel f\^ut bien notre sentiment en 
\'ecrivant ce travail que nous ne faisions que mettre par \'ecrit une succession de faits 
plus ou moins \emph{\'evidents}, 
ou tout au moins \emph{folkloriques}. Nous avons, de ce fait, \'et\'e surpris par les
cons\'equences de notre th\'eor\`eme principal, \`a savoir la version fonctorielle et 
multiplicative des isomorphismes HKR (voir \ref{cor}), isomorphismes certes 
bien connus (voir par exemple \cite{lo,sch,ye}), mais n'ayant pas
la r\'eputation de faits \emph{triviaux}, particuli\`erement pour les versions
globales valables sur des sch\'emas. La nouveaut\'e dans notre approche est 
de ne pas oublier que le complexe de Hochschild (d'une $k$-alg\`ebre commutative ou d'un 
$k$-sch\'ema) est muni de deux structures additionnelles, une multiplication et 
une action de $S^{1}$ (l'utilisation de la structure
multiplicative sur le complexe de Hochschild est 
clairement pr\'esente dans \cite{sch}). C'est l'existence de ces deux structures qui permet 
le lien naturel avec le th\'eorie de de Rham, et cela de mani\`ere essentiellement unique car 
ce lien est d\'eduit de propri\'et\'es universelles. \\

\bigskip

\textbf{Remerciements:} Nous remercions M. Hoyois qui a attir\'e notre attention 
sur le fait que la comparaison, annonc\'ee dans \cite{bena,tv',tv}, entre
fonctions sur les espaces de lacets d\'eriv\'es et homologie cyclique
\'etait probablement plus subtile que nous le pr\'etendions. Ce travail 
ne r\'epond pas directement \`a cette probl\'ematique. Il  montre cependant que l'on peut
directement obtenir une comparaison avec la th\'eorie de de Rham sans
passer par l'homologie cyclique, ce qui suffit pour les besoins de \cite{tv',tv} 
(et aussi visiblement de \cite{bena}). \\

\bigskip 

\textbf{Notations:} Tout au long de cet article $k$ d\'esigne un anneau 
commutatif \emph{de caract\'eristique nulle}. Les complexes de $k$-modules seront homologiquement indic\'es
et tous concentr\'es en degr\'es n\'egatifs par convention. La cat\'egorie
des complexes de $k$-modules, concentr\'es en degr\'es n\'egatifs, sera not\'ee
$C(k)$.  
Pour une $k$-dg-alg\`ebre $B$ (concentr\'ee en degr\'es n\'egatifs d'apr\`es nos conventions) 
on note $B-dg-mod$ la cat\'egorie des
$B$-dg-modules \`a gauche (de m\^eme, concentr\'es en degr\'es n\'egatifs). 
Cette cat\'egorie est munie d'une structure
de cat\'egorie de mod\`eles pour laquelle les \'equivalences
sont les quasi-isomorphismes et les fibrations sont les morphismes 
surjectifs en degr\'es strictement n\'egatifs. Nous noterons aussi 
$cdga$ la cat\'egorie des $k$-dg-alg\`ebres commutatives (toujours
en degr\'es n\'egatifs), que nous munirons de sa structure de mod\`eles
standard pour la quelle les \'equivalences
sont les quasi-isomorphismes et les fibrations sont les morphismes 
surjectifs en degr\'es strictement n\'egatifs (voir \cite{bogu,hi}). Nous utiliserons implicitement
que le foncteur de normalisation induit une \'equivalence entre
la cat\'egorie homotopique des $k$-alg\`ebres simpliciales commutatives et 
la cat\'egorie homotopique de $cdga$. Les techniques
de \cite{ss} impliquent que le foncteur de normalisation est alors
l'adjoint \`a gauche d'une \'equivalence de Quillen.

Nous noterons $S^{1}:=B\mathbb{Z}$ l'ensemble simplicial classifiant du groupe
ab\'elien $\mathbb{Z}$, que nous consid\'ererons toujours comme une 
groupe ab\'elien simplicial. En tant que groupe simplicial, $S^{1}$ peut op\'erer
sur tout objet dans une cat\'egorie simplicialement enrichie, et en particulier
dans une cat\'egorie de mod\`eles simpliciales $M$. Les objets $S^{1}$-\'equivariants
dans une telle cat\'egorie de mod\`eles $M$ forment une cat\'egorie
not\'es $S^{1}-M$. 

Nous utiliserons le langage, et des notions de bases, de la th\'eorie de d\'erivateurs
de Grothendieck (voir par exemple \cite{gr}). Le d\'erivateur associ\'e \`a une
cat\'egorie de mod\`eles $M$ sera not\'e $\mathbb{D}(M)$. Pour une
sous-cat\'egorie pleine $M_{0}$ d'une cat\'egorie de mod\`eles
$M$, stable par \'equivalences, nous noterons $\mathbb{D}(M_{0})$ le sous-d\'erivateur
plein de $\mathbb{D}(M)$ form\'e des objets de $M_{0}$. Toute adjonction 
de Quillen 
$$g : M \longrightarrow N \qquad M \longleftarrow N : f$$
induit une adjonction dans la 2-cat\'egorie des d\'erivateurs
$$\mathbb{L}g : \mathbb{D}(M)  \longrightarrow \mathbb{D}(N) \qquad \mathbb{D}(M) \longleftarrow
 \mathbb{D}(N) : \mathbb{R}f.$$
L'expression \emph{diagramme 2-commutatif de d\'erivateurs} fera r\'ef\'erence
\`a la donn\'e d'un diagramme de 1-morphismes munis de toutes les 2-isomorphismes
de coh\'erences n\'ecessaires. Ainsi, un carr\'e 2-commutatif est la donn\'ee non pas
de quatre 1-morphismes, mais bien de quatre 1-morphismes et un 2-isomorphisme
entre les deux compositions possibles.

Finalement, nous utiliserons aussi 
la notion de $\mathbb{S}$-cat\'egorie pour laquelle nous renvoyons \`a 
\cite{be} (et \`a \cite[\S 1]{tv} pour des propri\'et\'es plus avanc\'ees).

\section{Les $\epsilon$-dg-alg\`ebres}

On consid\`ere la $k$-dg-alg\`ebre $k[\epsilon]$, librement engendr\'ee
par un \'el\'ement $\epsilon$ en degr\'e $-1$ et avec la relation $\epsilon^{2}=0$. 
La $k$-alg\`ebre sous-jacente est $k[X]/X^{2}$, avec $deg(X)=-1$, et est munie
de la diff\'erentielle nulle. 

\begin{df}\label{d1}
La \emph{cat\'egorie des $\epsilon$-dg-modules} est la cat\'egorie
$k[\epsilon]-dg-mod$, des $k[\epsilon]$-dg-modules \`a gauche. Elle sera not\'ee
$\epsilon-dg-mod$. 
\end{df}

On remarque que $\epsilon-dg-mod$ n'est autre que la cat\'egorie des complexes
mixtes n\'egativement gradu\'es (au sens de \cite{lo} par exemple). \\

On munit $\epsilon-dg-mod$ de sa structure de mod\`eles usuelle o\`u les
\'equivalences sont les quasi-isomorphismes de complexes sous-jacents, et les fibrations
sont les morphismes surjectifs en degr\'es strictement n\'egatifs.  
La cat\'egorie $\epsilon-dg-mod$ est munie d'une structure mono\"\i dale sym\'etrique
induite par le produit tensoriel de complexes de $k$-modules. Plus pr\'ecis\'ement, 
pour $M$ et $N$ deux $\epsilon$-dg-modules on d\'efinit une structure
de $\epsilon$-dg-module sur le complexe $M\otimes_{k}N$ de la fa\c{c}on suivante.
Le complexe $M\otimes_{k}N$ est naturellement muni d'une structure de
$k[\epsilon]\otimes_{k} k[\epsilon]$-dg-module \`a gauche. On consid\`ere alors le morphisme
de $k$-dg-alg\`ebres
$$k[\epsilon] \longrightarrow k[\epsilon]\otimes_{k} k[\epsilon]$$
qui envoie $\epsilon$ sur $\epsilon\otimes 1 + 1\otimes \epsilon$. A travers ce morphisme
$M\otimes_{k} N$ est muni d'une structure de $\epsilon$-dg-module. On v\'erifie alors que les
contraintes d'associativit\'e, de sym\'etrie et d'unit\'e du produit tensoriel 
de complexes induisent des contraintes d'associativit\'e, de sym\'etrie et 
d'unit\'e pour la structure mono\"\i dale ainsi construite sur $\epsilon-dg-mod$. 
La cat\'egorie $\epsilon-dg-mod$ est ainsi munie d'une structure 
mono\"\i dale sym\'etrique que nous noterons simplement $\otimes$.

\begin{df}\label{d2}
La \emph{cat\'egorie des $\epsilon$-dg-alg\`ebres commutatives} est la cat\'egorie
des mono\"\i des associatifs, commutatifs et unitaires dans la cat\'egorie
monod\"\i dale $(\epsilon-dg-mod,\otimes)$. Elle sera not\'ee
$\epsilon-cdga$. 
\end{df}

En d'autres termes, un objet de $\epsilon-cdga$ consiste en une
$k$-dg-alg\`ebre commutative $A$, munie d'un morphisme 
de complexes de $k$-modules $\epsilon : A \longrightarrow A[1]$, tel que
pour tout $a$, $b$, \'el\'ements de $A$, de degr\'es respectifs $n$ et $m$, 
on ait
$$\epsilon(ab)=\epsilon(a)b + (-1)^{n}a\epsilon(b).$$

Nous allons maintenant munir $\epsilon-cdga$ d'une structure de cat\'egorie
de mod\`eles. Pour cela, nous
consid\'erons le foncteur d'oubli
$$\epsilon-cdga \longrightarrow C(k),$$
et nous d\'efinissons les \'equivalences (resp. les fibrations) dans 
$\epsilon-cdga$ comme les morphismes induisant des \'equivalences (resp. 
des fibrations) dans $C(k)$. Ainsi, une \'equivalence 
dans $\epsilon-cdga$ est un morphisme induisant un quasi-isomorphisme 
sur le complexes sous-jacent. De m\^eme, une fibration dans $\epsilon-cdga$ est un morphisme 
qui est surjectif en tout degr\'e strictement n\'egatif.\\
Ce foncteur d'oubli poss\`ede un adjoint
\`a gauche. En effet, il est facile de voir que le foncteur d'oubli commute \`a tout type de limites
ainsi qu'aux colimites filtrantes. Comme les cat\'egories $\epsilon-cdga$ et $cdga$
sont des cat\'egories localement pr\'esentables, l'existence d'un adjoint \`a gauche
est assur\'ee par le th\'eor\`eme d'existence de Freyd. 
Nous noterons 
$$\epsilon : cdga \longrightarrow \epsilon-cdga$$
l'adjoint \`a gauche du foncteur d'oubli, dont nous allons maintenant donner une
construction plus explicite. Soit $A\in cdga$, et notons $\Omega_{A}^{1}$
le $A$-dg-module corepr\'esentant le foncteur des d\'erivations (au sens dg). 
Ce $A$-dg-module est librement engendr\'e sur $A$ par des symb\^oles
$\partial(a)$, avec $a\in A$, $deg(\partial(a))=deg(a)$, et avec les relations
$$\partial(ab)=(-1)^{deg(a).deg(b)}b\partial(a) + (-1)^{deg(a)}a\partial(b).$$
On note $DR(A)$ la $A$-dg-alg\`ebre commutative libre sur $\Omega_{A}^{1}[1]$
$$DR(A):=\oplus_{n}(\Omega_{A}^{1}[1])^{\otimes_{A}\, n}/\Sigma_{n}.$$
On munit enfin cette $cdga$ d'une $\epsilon$-structure en d\'ecr\'etant que
$$\epsilon : A \subset DR(A) \longrightarrow \Omega_{A}^{1} \subset DR(A)[-1]$$
est la d\'erivation universelle $a \mapsto \partial(a)$, et en prolongeant par 
multiplicativit\'e. La $k$-dg-alg\`ebre commutative, munie de 
$\epsilon$, est un objet de $\epsilon-cdga$ not\'e $\epsilon(A)$. De plus, 
la $cdga$ sous-jacente \`a $\epsilon(A)$ est $DR(A)$, et le morphisme naturel
$A \longrightarrow DR(A)$ induit une bijection
$$Hom_{\epsilon-cdga}(DR(A),B) \longrightarrow Hom_{cdga}(DR(A),B) \longrightarrow 
Hom_{cdga}(A,B).$$
Ainsi, $A \mapsto \epsilon(A)$ est un adjoint \`a gauche du foncteur
d'oubli.

\begin{prop}\label{p1}
Les notions ci-dessus de fibrations et d'\'equivalences munissent 
$\epsilon-cdga$ d'une structure de cat\'egorie de mod\`eles. De plus, le 
foncteur d'oubli $\epsilon-cdga \longrightarrow cdga$ est de Quillen 
\`a droite.
\end{prop}

\textit{Preuve:} Il s'agit de relever la structure de mod\`eles 
sur $cdga$ le long du foncteur d'oubli $\epsilon-cdga \longrightarrow cdga$.\\ 
Notons $I_{0}$ et $J_{0}$ des ensembles g\'en\'erateurs de cofibrations et cofibrations triviales
dans $cdga$. On d\'efinit $I:=\epsilon(I_{0})$, l'image de $I_{0}$ par le foncteur
$\epsilon$. De m\^eme, on pose $J:=\epsilon(J_{0})$. On applique alors
le th\'eor\`eme 2.1.19 de \cite{ho}. Comme le foncteur d'oubli $\epsilon-cdga \longrightarrow cdga$ 
refl\`ete les limites, les colimites, les fibrations et les \'equivalences, on voit que
pour v\'erifier les conditions de ce th\'eor\`eme il suffit de montrer que
$J \subset W$. D'apr\`es la construction explicite du foncteur $\epsilon$ que nous
avons donn\'e pr\'ec\'edemment on voit qu'il suffit de montrer que
$A \mapsto \Omega_{A}^{1}$ (en tant que foncteur 
$cdga \longrightarrow C(k)$) transforme cofibrations triviales en quasi-isomorphismes.
Pour cela il suffit de montrer que pour $A \longrightarrow B$ une cofibration triviale
de $cdga$, le morphisme induit
$$\Omega_{A}^{1}\otimes_{A}B \longrightarrow \Omega_{B}^{1}$$
est une cofibration triviale de $B$-dg-modules. Pour cela, donnons nous  un 
diagramme commutatif de $B$-dg-modules
$$\xymatrix{
\Omega_{A}^{1}\otimes_{A}B \ar[r] \ar[d] & M \ar[d] \\
\Omega_{B}^{1} \ar[r] & N,}$$
avec $M \longrightarrow N$ une fibration. Par propri\'et\'e universelle des
dg-modules $\Omega^{1}$, on voit que ce diagramme correspond \`a un diagramme
commutatif dans $cdga/B$
$$\xymatrix{
A \ar[r] \ar[d] & B\oplus M \ar[d] \\
B \ar[r] & B \oplus N,}$$
o\`u $B\oplus E$ est l'extension de carr\'e nulle triviale de $B$ par 
le dg-module $E$. Comme $A \longrightarrow B$ est une cofibration triviale et que
$B\oplus M \longrightarrow B \oplus N$ est une fibration, 
il existe $B \longrightarrow B\oplus M$ un rel\`evement
dans $cdga/B$. Par adjonction on voit que cela implique l'existence
d'un rel\`evement $\Omega_{B}^{1} \longrightarrow M$
de $B$-dg-modules. Ceci montre donc que 
$\Omega_{A}^{1}\otimes_{A}B \longrightarrow \Omega_{B}^{1}$ rel\`eve 
les fibrations et donc est une cofibration triviale. En particulier, 
le morphisme $\Omega_{A}^{1} \longrightarrow \Omega_{B}^{1}$
est une \'equivalence, ce qu'il nous fallait montrer.
\hfill $\Box$ \\

Le foncteur $\epsilon$, retreint 
\`a la sous-cat\'egorie de $cdga$ form\'ee des $k$-alg\`ebres non-dg, poss\`ede
l'interpr\'etation plaisante suivante. Soit $A$ une $k$-alg\`ebre
commutative et 
$$DR(A)=Sym_{A}(\Omega_{A}^{1}[1])\simeq \oplus_{n} \Omega_{A}^{n}[n]$$ 
sa dg-alg\`ebre de de Rham. La diff\'erentielle de de Rham 
munit $DR(A)$ d'une structure de $\epsilon-cdga$ qui n'est autre que 
$\epsilon(A)$. En d'autres termes, $DR(A)$ muni de sa diff\'erentielle de
de Rham est la $\epsilon$-dg-alg\`ebre commutative libre engendr\'ee par $A$. 
De plus, si $A$ est cofibrante en tant qu'objet de $cdga$ (e.g. $A$ est une
$k$-alg\`ebre de polyn\^omes) alors $\epsilon(A)$ est cofibrante
dans $\epsilon-cdga$. 

\begin{prop}\label{p2}
Notons 
$$\mathbb{L}\epsilon : Ho(cdga) \longrightarrow Ho(\epsilon-cdga)$$
le foncteur d\'eriv\'e \`a gauche de $\epsilon$. Si 
$A$ est une $k$-alg\`ebre commutative lisse alors le morphisme naturel
$$\mathbb{L}\epsilon(A) \longrightarrow \epsilon(A)$$
est un isomorphisme dans $Ho(\epsilon-cdga)$.
\end{prop}

\textit{Preuve:} Soit $QA \longrightarrow A$ un mod\`ele
cofibrant pour $A$ dans $cdga$. Le morphisme en question est repr\'esent\'e
par
$$Sym_{QA}(\Omega_{QA}^{1}[1]) \longrightarrow Sym_{A}(\Omega_{A}^{1}[1]).$$
Ainsi, ce morphisme est un quasi-isomorphisme si et seulement si le morphisme induit
$$\Omega_{QA}^{1} \longrightarrow \Omega_{A}^{1}$$
est un quasi-isomorphisme de complexes.
Or, $\Omega_{QA}^{1}$ est un mod\`ele pour le complexe cotangent 
$\mathbb{L}_{A}$ de $A$, comme cela se voit en utilisant 
l'\'equivalence de Quillen entre $cdga$ et $k$-alg\`ebres simpliciales commutatives, ainsi
que la caract\'erisation du complexe cotangent en termes de d\'erivations (voir par exemple \cite{hagII}).
Le morphisme ci-dessus est alors isomorphe, dans $Ho(C(k))$, au 
morphisme naturel
$$\mathbb{L}_{A} \longrightarrow \Omega^{1}_{A}.$$
Or, comme $A$ est lisse ce morphisme est bien un quasi-isomorphisme.
\hfill $\Box$ \\

\section{Alg\`ebres simpliciales $S^{1}$-\'equivariantes et $\epsilon$-dg-alg\`ebres}

Notons $sk-CAlg$ la cat\'egorie des $k$-alg\`ebres commutatives simpliciales, et consid\'erons
$S^{1}-sk-CAlg$ la cat\'egorie des objets de $sk-CAlg$ avec une action du 
groupe simplicial $S^{1}$. La cat\'egorie $sk-CAlg$ est munie d'une structure de mod\`eles
pour la quelle les fibrations et les \'equivalences sont d\'efinies sur les
ensembles simpliciaux sous-jacents. Comme $sk-CAlg$ est une cat\'egorie
de mod\`eles simpliciale et engendr\'ee par cofibration la cat\'egorie
$S^{1}-sk-CAlg$ est elle-m\^eme munie d'une structure projective o\`u les fibrations et 
\'equivalences sont d\'efinis dans $sk-CAlg$. 

Dans cette section nous allons construire une \'equivalence de d\'erivateurs
$$\phi : \mathbb{D}(S^{1}-sk-CAlg) \longrightarrow \mathbb{D}(\epsilon-dg-alg)$$
ainsi qu'un $2$-isomorphisme $h$ faisant commuter le diagramme suivant
(en tant que diagramme dans la $2$-cat\'egorie des d\'erivateurs)

$$\xymatrix{
\mathbb{D}(S^{1}-sk-CAlg) \ar[d] \ar[r]^-{\phi}  & \ar[d] \mathbb{D}(\epsilon-dg-mod) \\
\mathbb{D}(sk-CAlg) \ar[r]_-{N} & \mathbb{D}(cdga).}$$
Dans ce diagramme les morphismes verticaux sont induits par les foncteurs
d'oubli 
$$S^{1}-sk-CAlg \longrightarrow sk-CAlg \qquad \epsilon-cdga \longrightarrow cdga,$$
et le foncteur $N$ est le foncteur de normalisation, adjoint \`a droite d'une \'equivalence de Quillen
(voir \cite{ss})
$$N : sk-CAlg \longrightarrow cdga.$$

\textbf{Premi\`ere \'etape:} 
On consid\`ere le groupe simplicial $BS^{1}$ comme une
$\mathbb{S}$-cat\'egorie avec une unique objet $*$, dont le
mono\"\i de des endomorphismes est $S^{1}$. On choisit alors une cat\'egorie
$\mathcal{C}$, et un diagramme de $\mathbb{S}$-cat\'egories
$$\xymatrix{\mathcal{C} \ar[r]^-{i} & T & \ar[l]_{j} BS^{1},}$$
tel que, d'une part $j$ soit une \'equivalence de $\mathbb{S}$-cat\'egories, 
et $i$ fasse de $T$ une localisation de $\mathcal{C}$ le long de tous
ses morphismes (au sens par exemple de \cite[\S 1.2]{tv}). On pourra, par exemple, prendre
pour $\mathcal{C}$ la cat\'egorie cyclique $\Lambda$.
On sait alors qu'il existe une chaine d'adjonctions
de Quillen 
$$sk-CAlg^{\mathbb{C}} \longleftrightarrow sk-CAlg^{T} \longleftrightarrow sk-CAlg^{BS^{1}}=S^{1}-sk-CAlg.$$
Ces adjunctions induisent des \'equivalences de d\'erivateurs (voir par exemple
\cite[\S 2.3.2]{tv0} ou encore \cite[\S 1.2]{tv})
$$\xymatrix{
\mathbb{D}(S^{1}-sk-CAlg) \ar[r]^-{j^{*}} & \mathbb{D}(sk-CAlg^{T}) & \ar[l]_-{i^{*}} \mathbb{D}_{loc}(sk-CAlg^{\mathbb{C}}),}$$
o\`u $\mathbb{D}_{loc}(sk-CAlg^{\mathbb{C}})$ d\'esigne le sous-d\'erivateur plein de 
$\mathbb{D}(sk-CAlg^{\mathbb{C}})$ form\'e des diagrammes $\mathcal{C} \longrightarrow sk-CAlg$ 
qui envoie tous les morphismes de $\mathcal{C}$ sur des \'equivalences. Fixons-nous
un objet $x\in \mathcal{C}$, alors le foncteurs ci-dessus commutent clairement 
\`a l'\'evaluation en $x$ et au point de base de $BS^{1}$, et on dispose ainsi 
de deux carr\'es $2$-commutatifs
$$\xymatrix{
\mathbb{D}(S^{1}-sk-CAlg) \ar[dr] \ar[r]^-{j^{*}} & \mathbb{D}(sk-CAlg^{T}) \ar[d]^-{ev_{i(x)}} & \ar[l]_-{i^{*}} \mathbb{D}_{loc}(sk-CAlg^{\mathbb{C}}) \ar[dl]^-{ev_{x}} \\
 & \mathbb{D}(SEns), & }
$$
o\`u les morphismes horizontaux sont des \'equivalences, et le morphisme vertical de gauche est le
foncteur d'oubli. On construit ainsi un diagramme 2-commutatif de d\'erivateurs
$$\xymatrix{
\mathbb{D}(S^{1}-sk-CAlg) \ar[rr]^-{\phi_{1}}  \ar[dr] & & \mathbb{D}_{loc}(sk-CAlg^{\mathbb{C}}) \ar[dl]^-{ev_{x}}
\\
 & \mathbb{D}(SEns). & }$$
Comme $BS^{1}$ est simplement connexe, il n'est pas difficile de v\'erifier que ce diagramme ne d\'epend pas,
\`a \'equivalence pr\`es, du choix du point $x\in \mathcal{C}$.  \\

\textbf{Seconde \'etape:} Consid\'erons le foncteur de normalisation
$$N : sk-CAlg \longrightarrow cdga$$
que l'on sait \^etre l'adjoint \`a droite d'une \'equivalence de Quillen (voir \cite{ss}).
Il induit donc un nouvel adjoint \`a droite d'une \'equivalence Quillen
$$N : sk-CAlg^{\mathcal{C}} \longrightarrow cdga^{\mathcal{C}},$$
qui induit, \`a sont tour, une \'equivalence de d\'erivateurs
$$\mathbb{D}_{loc}(sk-CAlg^{\mathcal{C}}) \longrightarrow \mathbb{D}_{loc}(cdga^{\mathcal{C}}).$$
Cette \'equivalence vient avec un 2-isomorphisme naturel faisant commutater le diagramme
suivant
$$\xymatrix{
\mathbb{D}_{loc}(sk-CAlg^{\mathcal{C}}) \ar[rr]^{\phi_{2}} \ar[dr]_-{ev_{x}} & & \mathbb{D}_{loc}(cdga^{\mathcal{C}})
\ar[dl]^-{ev_{x}} \\
 & \mathbb{D}(cdga). & }$$

\textbf{Troisi\`eme \'etape:} Dans cette \'etape, et la suite,  nous aurons besoin de travailler
momentan\'ement avec des complexes non born\'es. Lorsque cela sera le cas nous l'indiquerons
par un indice $(-)_{\infty}$. Ainsi, $C(k)_{\infty}$, $cdga_{\infty}$ \dots d\'esignera la
cat\'egorie des complexes non born\'es, des $k$-dg-al\`ebres commutatives non born\'ees \dots.

L'inclusion des complexes en degr\'es n\'egatifs dans les complexes non born\'es induit un 
morphisme de d\'erivateurs 
$$\mathbb{D}_{loc}(cdga^{\mathcal{C}}) \longrightarrow \mathbb{D}_{loc}(cdga_{\infty}^{\mathcal{C}}).$$
Ce morphisme est pleinement fid\`ele et identifie le membre de gauche au sous-d\'erivateur
des objets cohomologiquement concentr\'es en degr\'es n\'egatifs. 

On consid\`ere le functor des sections globales
$$\Gamma : cdga_{\infty}^{\mathcal{C}} \longrightarrow cdga_{\infty}.$$
Ce foncteur est de Quillen \`a droite lorsque l'on munit 
$cdga_{\infty}^{\mathcal{C}}$ de sa structure injective, pour laquelle les
cofibrations et les \'equivalences sont d\'efinies termes \`a termes. Pour tout
$A \in cdga_{\infty}^{\mathcal{C}}$, le morphisme unit\'e $\underline{k} \longrightarrow A$, 
o\`u $\underline{k}$ est le diagramme constant, induit un morphisme dans $cdga$
$$\Gamma(\underline{k}) \longrightarrow \Gamma(A).$$
Ainsi, si $R$ d\'esigne un foncteur de remplacement fibrant dans $cdga_{\infty}^{\mathcal{C}}$, 
on dispose d'un morphisme
$$\Gamma(R(\underline{k})) \longrightarrow \Gamma(R(A)).$$
Nous noterons $B:=\Gamma(R(\underline{k})) \in cgda$. La construction $A \mapsto \Gamma(R(A))$
d\'efinit ainsi un morphisme de d\'erivateurs
$$\mathbb{D}_{loc}(cdga_{\infty}^{\mathcal{C}}) \longrightarrow
\mathbb{D}(B-cdga),$$
o\`u $B-cdga$ d\'esigne la cat\'egorie comma $B/cdga$. 

Rappelons que, lors de la seconde \'etape, nous nous sommes fix\'es un diagramme
de $\mathbb{S}$-cat\'egories
$$\xymatrix{\mathcal{C} \ar[r]^-{i} & T & \ar[l]_-{j} BS^{1}.}$$
Ce diagramme induit des isomorphismes de $k$-alg\`ebres gradu\'ees commutatives de cohomologie
$$H^{*}(B)=H^{*}(\mathbb{R}\Gamma(\underline{k}))\simeq H^{*}(\mathcal{C},k)
\simeq H^{*}(BS^{1},k)\simeq k[u],$$
o\`u $deg(u)$ et correspond au g\'en\'erateur de $H^{2}(K(\mathbb{Z},2),k)$ 
donn\'e par l'inclusion standard $\mathbb{Z} \subset k$. Le choix 
d'un $2$-cocycle $u' \in Z^{2}(B)$ qui est un repr\'esentant de $u$ d\'etermine
un quasi-isomorphisme de cdga $k[u]\longrightarrow B$. Ce quasi-isomorphisme, 
consid\'er\'e \`a homotopie pr\`es, ne d\'epend pas du choix
de $u'$. Il induit ainsi une \'equivalence de Quillen 
$$B-cdga_{\infty} \longrightarrow k[u]-cdga_{\infty}$$
dont le morphisme correspond de d\'erivateurs
$$\mathbb{D}(B-cdga_{\infty}) \longrightarrow \mathbb{D}(k[u]-cdga_{\infty})$$
est une \'equivalence, d\'etermin\'ee \`a 2-isomorphisme unique pr\`es.

Nous avons ainsi construit un morphisme de d\'erivateurs
$$\phi_{3} : \mathbb{D}_{loc}(cdga^{\mathcal{C}}) \hookrightarrow \mathbb{D}(cdga_{\infty}^{\mathcal{C}})
\longrightarrow \mathbb{D}(B-cdga_{\infty}) \longrightarrow \mathbb{D}(k[u]-cdga_{\infty}).$$
Ce morphisme entre dans un diagramme 
$$\xymatrix{
\mathbb{D}_{loc}(cdga^{\mathcal{C}}) \ar[r]^-{\phi_{3}} \ar[d]_-{q} & \mathbb{D}(k[u]-cdga_{\infty}) \ar[d]^-{p} \\
\mathbb{D}(cdga) \ar[r]_-{i} & \mathbb{D}(cdga^{\infty}),}$$
qui n'est pas 2-commutatif et demande quelques explications suppl\'ementaires.
Le morphisme $q$ est induit par le foncteur d'oubli, et $i$ par l'inclusion
naturelle. Cependant, $p$ n'est pas le foncteur d'oubli. Il est 
induit par le foncteur de Quillen \`a gauche 
$$k[u]-cdga_{\infty} \longrightarrow cdga_{\infty}$$
qui envoie $A'$ sur $k\otimes_{k[u]}A'$. Par adjonction, il n'est pas difficile 
de voir qu'il existe un 2-morphisme
$$u : p\circ \phi_{3} \Rightarrow i\circ q,$$
qui n'est pas un 2-isomorphisme en g\'en\'eral. Il le deviendra lorsque les foncteurs
seront restreint \`a certains sous-d\'erivateurs d'objets born\'es, comme nous
allons maintenant le voir.

Nous noterons $\mathbb{D}^{+}_{loc}(cdga^{\mathcal{C}})$ (resp.
$\mathbb{D}^{+}(k[u]-cdga_{\infty})$) le sous-d\'erivateur form\'ee
des objets dont les complexes sous-jacents sont cohomologiquement born\'es
\`a gauche (i.e. $H^{i}$ s'annule pour tout $i$ suffisament petit). Nous remarquerons ici que la
restriction de $\phi_{3}$
$$\phi_{3}^{+} : \mathbb{D}^{+}_{loc}(cdga^{\mathcal{C}}) \longrightarrow \mathbb{D}^{+}(k[u]-cdga_{\infty})$$
est pleinement fid\`ele.  De plus, la restriction du 2-morphisme $u$ ci-dessus induit un 
2-isomorphisme
$$u^{+} : p\circ \phi^{+}_{3} \Rightarrow i\circ q.$$
Pour voir cela, il nous faut revenir \`a l'adjonction de Quillen
$$cdga_{\infty} \longleftrightarrow cdga^{\mathcal{C}}_{\infty}.$$ 
Notons $A:=R(\underline{k})$ un mod\`ele fibrant de $\underline{k}$ dans $cdga^{\mathcal{C}}_{\infty}$, et
$B=\Gamma(A)$. 
On consid\`ere l'adjonction induite
$$f : B-cdga_{\infty} \longleftrightarrow A/cdga_{\infty}^{\mathcal{C}} : \Gamma.$$
Le foncteur $f$ envoie un objet $B' \in B-cdga_{\infty}$ sur 
$A\otimes_{B}B'$, o\`u $A$ est consid\'er\'ee comme une $B$-dg-alg\`ebre commutative
\`a travers le morphisme $B \longrightarrow A$, adjoint de l'identit\'e
$B=\Gamma(A)$ (nous identifions ici les objets de
$B-cdga_{\infty}$ avec leur diagrammes constants correspondants). Il nous faut montrer 
que pour tout objet $A' \in A/cdga_{\infty}^{\mathcal{C}}$, cohomologiquement 
born\'e \`a gauche, le morphisme d'adjonction
$$A\otimes_{B}^{\mathbb{L}}\mathbb{R}\Gamma(A') \longrightarrow  A'$$
est un isomorphisme dans $Ho(cdga_{\infty}^{\mathcal{C}})$. Pour cela on peut oublier
la structure d'alg\`ebres et consid\'er\'e que $A'$ est un $A$-dg-module. Dans ce cas, un 
argument de type d\'ecomposition de Postnikov sur $A'$ ram\`ene le probl\`eme
au cas o\`u $A'$ est un diagramme constant associ\'e \`a un $k$-modules $M$. On a alors
$$H^{*}(\mathbb{R}\Gamma(M))\simeq H^{*}(\mathcal{C},M)\simeq H^{*}(K(\mathbb{Z},2),M)\simeq k[u]\otimes_{k} M.$$ 
On a donc
$$A\otimes_{B}^{\mathbb{L}}\mathbb{R}\Gamma(M)\simeq A\otimes_{k[u]}^{\mathbb{L}}k[u]\otimes_{k}M \simeq M.$$
En revenant aux d\'efinitions de nos foncteurs, ceci montre aussi que $u^{+}$ est un 2-isomorphisme
(nous laissons la v\'erification au lecteur).
On peut de plus caract\'eriser l'image de $\phi_{3}^{+}$. En effet, 
l'objet $k[u]$ engendre une t-structure sur le d\'erivateur $\mathbb{D}(k[u]-dg-mod_{\infty})$, 
dont la partie n\'egative est engendr\'ee par $k[u]$ par colimites homotopiques. 
L'image essentielle de $\phi_{3}^{+}$ consiste alors en tous les objets
de $\mathbb{D}(k[u]-cdga_{\infty})$ dont l'objet sous-jacent dans $\mathbb{D}(k[u]-dg-mod_{\infty})$
est d'amplitude $[n,0]$ pour un certain entier $n$. \\

\textbf{Quatri\`eme \'etape:} De mani\`ere analogue \`a l'\'etape pr\'ec\'edente, nous
construisons un diagramme, qui n'est pas 2-commutatif, de d\'erivateurs
$$\xymatrix{
\mathbb{D}(\epsilon-cdga) \ar[r]^-{\phi_{4}} \ar[d]_-{r} & \mathbb{D}(k[u]-cdga_{\infty}) \ar[d]^-{p} \\
\mathbb{D}(cdga) \ar[r]_-{i} & \mathbb{D}(cdga^{\infty}),}$$
tel que la restriction de $\phi_{4}$ \`a $\mathbb{D}^{+}\epsilon-cdga)$, le sous-d\'erivateurs
des objets cohomologiquement born\'es, soit pleinement fid\`ele. Nous construisons aussi
un 2-morphisme $v : p\circ \phi_{4} \Rightarrow i\circ r$
qui induira un 2-isomorphisme par restriction
$$v^{+} : p\circ \phi^{+}_{4} \Rightarrow i\circ r.$$
La construction de 
$\phi_{4}$ et de $v$ est tout \`a fait analogue \`a celle de $\phi_{3}$ et de $u$, nous
nous contenterons donc de l'esquisser. Nous consid\'ererons
$\epsilon-dg-mod$, muni de sa structure \emph{injective}, pour laquelle les cofibrations
et les \'equivalences sont d\'efinies dans $C(k)$. Cette structure de mod\`eles 
reste une structure de mod\`eles mono\"\i dales au sens de \cite{ho}, et 
elle induit une structure injective sur $\epsilon-cdga$, pour laquelle les cofibrations et 
les \'equivalences sont d\'efinies dans $cdga$. On consid\`ere alors le foncteur
$$\Gamma_{\epsilon} : \epsilon-cdga \longrightarrow cdga_{\infty},$$
qui envoie $A \in \epsilon-cdga$ sur $\underline{Hom}_{\epsilon-dg-mod_{\infty}}(k,A)$, o\`u 
$\underline{Hom}_{\epsilon-dg-mod_{\infty}}$ d\'esigne le complexe, non born\'e, des morphismes
de $\epsilon$-dg-modules. Le foncteur $\Gamma_{\epsilon}$ est de Quillen \`a droite, et la construction
$A \mapsto \Gamma_{\epsilon}(R(A))$, o\`u $R$ est un remplacement fibrant dans $\epsilon-cdga$, fournit
un morphisme de d\'erivateurs
$$\mathbb{D}(\epsilon-cdga) \longrightarrow \mathbb{D}(B'-cdga_{\infty}),$$
avec $B':=\Gamma_{\epsilon}(R(k))$. Or, $H^{*}(B')\simeq Ext^{*}_{k[\epsilon]}(k,k)\simeq k[u]$. Ainsi, 
il existe un quasi-isomorphisme de dg-alg\`eres commutatives $k[u] \longrightarrow B'$, bien 
d\'etermin\'e \`a homotopie pr\`es. Ce quasi-isomorphisme fournit une \'equivalence
de d\'erivateurs $\mathbb{D}(B'-cdga_{\infty}) \longrightarrow \mathbb{D}(k[u]-cdga_{\infty})$.
Le morphisme $phi_{4}$ est par d\'efinition l'\'equivalence compos\'ee
$$\xymatrix{
\mathbb{D}(\epsilon-cdga) \ar[r] & \mathbb{D}(B'-cdga_{\infty}) \ar[r]^-{\sim} & 
\mathbb{D}(k[u]-cdga_{\infty}).}$$
Nous laissons le soins au lecteur de v\'erifier que $\phi_{4}^{+}$ est bien 
pleinement fid\`ele, et que son image essentielle co\"\i ncide avec celle
de $\phi_{3}^{+}$ (les arguments sont les m\^emes que pour $\phi_{3}$). \\

\textbf{Cinqui\`eme et derni\`ere \'etape:} D'apr\`es les deux derni\`eres \'etapes 
nous avons construit un diagramme 2-commutatif de d\'erivateurs
$$\xymatrix{
\mathbb{D}_{loc}^{+}(cdga^{\mathcal{C}}) \ar[r]^-{\phi_{3}^{+}} \ar[d] & \mathbb{D}(k[u]-cdga_{\infty}) \ar[d] 
& \ar[l]_-{\phi_{4}^{+}} \mathbb{D}^{+}(\epsilon-cdga) \ar[d] \\
\mathbb{D}^{+}(cdga) \ar[r] & \mathbb{D}(cdga^{\infty}) & \ar[l] \mathbb{D}^{+}(cdga).}$$
Les morphismes $\phi_{3}^{+}$ et $\phi_{4}^{+}$ \'etant pleinement fid\`eles avec la
m\^eme image essentielle on trouve un nouveau diagramme 2-commutatif
$$\xymatrix{
\mathbb{D}_{loc}^{+}(cdga^{\mathcal{C}}) \ar[r]^-{\phi_{34}^{+}} \ar[d] & \mathbb{D}^{+}(\epsilon-cdga) \ar[d] \\
\mathbb{D}^{+}(cdga) \ar[r]_-{id} & \mathbb{D}^{+}(cdga).}$$
Il n'est pas difficile de remarquer que ce diagramme se compl\`ete, de mani\`ere unique (\`a isomorphisme unique
pr\`es), en un diagramme 2-commutatif
$$\xymatrix{
\mathbb{D}_{loc}(cdga^{\mathcal{C}}) \ar[r]^-{\phi_{34}} & \mathbb{D}(\epsilon-cdga) \\
\mathbb{D}_{loc}^{+}(cdga^{\mathcal{C}})\ar[u] \ar[r]^-{\phi_{34}^{+}} \ar[d] & 
\mathbb{D}^{+}(\epsilon-cdga) \ar[d] \ar[u]\\
\mathbb{D}^{+}(cdga) \ar[r]_-{id} & \mathbb{D}^{+}(cdga),}$$
avec $\phi_{34}$ une \'equivalence de d\'erivateurs. En effet, il suffit pour cela de remarquer
que dans $\mathbb{D}_{loc}(cdga^{\mathcal{C}})$ et $\mathbb{D}(\epsilon-cdga)$, tout objet 
est limite de sa tour de Postnikov. Ainsi, $\phi_{34}$ est d\'efini est prenant l'image par $\phi_{34}^{+}$
des tours de Postnikov de $\mathbb{D}_{loc}(cdga^{\mathcal{C}})$, 
puis en en prenant leurs limites dans le d\'erivateur $\mathbb{D}(\epsilon-cdga)$. Nous laissons
les d\'etails au lecteur de cette construction formelle. \\

La conclusion de ces cinq \'etapes et la construction d'un diagramme 2-commutatif
de d\'erivateurs
$$\xymatrix{
\mathbb{D}(S^{1}-sk-CAlg) \ar[r]^-{\phi} \ar[d] & \mathbb{D}(\epsilon-cdga) \ar[d] \\
\mathbb{D}(sk-CAlg) \ar[r]_-{N} & \mathbb{D}(cdga),}$$
o\`u $\phi:=\phi_{34}\circ \phi_{2}\circ \phi_{1}$, le morphisme $N$ est induit par
le foncteur de normalisation, et les morphismes verticaux sont les
oublis.

\section{Le th\'eor\`eme de comparaison et quelques applications}

Nous sommes maintenant en mesure de d'\'enoncer et de d\'emontrer notre
r\'esultat principal. 

\begin{thm}\label{t1}
\begin{enumerate}
\item Les deux morphismes de d\'erivateurs
$$\mathbb{D}(sk-CAlg) \longrightarrow \mathbb{D}(\epsilon-cdga),$$
qui envoient respectivement $A$ sur $\phi(S^{1}\otimes^{\mathbb{L}}_{k}A)$ et 
sur $\mathbb{L}\epsilon(N(A))$, sont naturellement isomorphes.
\item Soit $\mathbb{D}(k-CAlg^{sm})$ le sous-d\'erivateur plein de
$\mathbb{D}(sk-CAlg)$ form\'e des $k$-alg\`ebres commutatives lisses et sur $k$.
Alors les deux morphismes de d\'erivateurs
$$\mathbb{D}^{sm}(k-CAlg) \longrightarrow \mathbb{D}(\epsilon-cdga),$$
qui envoient respectivement $A$ sur $\phi(S^{1}\otimes_{k}A)$ et 
sur $\epsilon(A)$, sont naturellement isomorphes.
\end{enumerate}
\end{thm}

\textit{Preuve:} Tout d'abord, $(2)$ d\'ecoule de $(1)$ et du fait que pour 
$A$ lisse sur $k$ les deux morphismes
$$\mathbb{L}\epsilon(A) \longrightarrow \epsilon(A) \qquad
S^{1}\otimes^{\mathbb{L}}_{k}A \longrightarrow S^{1}\otimes_{k}A$$
sont des \'equivalences (\`a cause de la proposition \ref{p2}, et parceque 
$A$ est plate sur $k$). Pour $(1)$, on revient au diagramme 2-commutatif
$$\xymatrix{
\mathbb{D}(S^{1}-sk-CAlg) \ar[r]^-{\phi} \ar[d] & \mathbb{D}(\epsilon-cdga) \ar[d] \\
\mathbb{D}(sk-CAlg) \ar[r]_-{N} & \mathbb{D}(cdga),}$$
et on remarque que les adjoints \`a gauche des oublis sont respectivement 
induits par les foncteurs de Quillen \`a gauche
$$\epsilon : cdga \longrightarrow \epsilon-cdga \qquad
S^{1}\otimes_{k} - : sk-CAlg \longrightarrow S^{1}-sk-CAlg.$$
Comme les morphismes $\phi$ et $N$ sont des \'equivalences le diagramme 
induit
$$\xymatrix{
\mathbb{D}(S^{1}-sk-CAlg) \ar[r]^-{\phi}  & \mathbb{D}(\epsilon-cdga)  \\
\mathbb{D}(sk-CAlg) \ar[u]^-{S^{1}\otimes^{\mathbb{L}} -} \ar[r]_-{N} & \mathbb{D}(cdga) 
\ar[u]_-{\mathbb{L}\epsilon},}$$
est naturellement 2-commutatif. \hfill $\Box$ \\

Pour le corollaire suivant, rappelons que pour $X$ un 
$k$-sch\'ema s\'epar\'e (\`a diagonale affine suffirait), 
on dispose de son sch\'ema d\'eriv\'e
$LX:=\mathbb{R}Map(S^{1},X)$ (voir \cite{tv,tv'} ainsi que 
\cite[4.3.1]{to}). Il vient avec 
une projection naturelle $LX \longrightarrow X$, qui fait
de $LX$ le spectre relatif du faisceau de $k$-alg\`ebres simpliciales
commutatives $S^{1}\otimes^{\mathbb{L}}\mathcal{O}_{X}$ sur $X$. Le groupe
simplicial $S^{1}$ op\`ere naturellement sur $LX$ en agissant sur lui-m\^eme
par translations. 

\begin{cor}\label{cor}
Soit $X$ un $k$-sch\'ema de type fini et s\'epar\'e sur $Spec\, k$.
\begin{enumerate}
\item Il existe un isomorphisme dans la cat\'egorie homotopique
des faisceaux de $\epsilon-cdga$ sur $X$
$$\phi(S^{1}\otimes^{\mathbb{L}}\mathcal{O}_{X})\simeq \mathbb{L}\epsilon(\mathcal{O}_{X}).$$
\item Il existe un isomorphisme dans la cat\'egorie homotopique
des faisceaux de $\mathcal{O}_{X}-cdga$ sur $X$
$$\mathcal{O}_{X}\otimes^{\mathbb{L}}_{\mathcal{O}_{X}\otimes_{k}^{\mathbb{L}}
\mathcal{O}_{X}}\mathcal{O}_{X} \simeq Sym_{\mathcal{O}_{X}}(\mathbb{L}_{X/k}[1]),$$
o\`u $\mathbb{L}_{X/k}$ est le complexe cotangent de $X$ relativement \`a $k$ au sens
de \cite{il}.
\item Si $X$ est lisse sur $k$, alors il existe 
isomorphisme dans la cat\'egorie homotopique
des faisceaux de $\mathcal{O}_{X}-cdga$ sur $X$
$$\mathcal{O}_{X}\otimes^{\mathbb{L}}_{\mathcal{O}_{X}\otimes_{k}^{\mathbb{L}}
\mathcal{O}_{X}}\mathcal{O}_{X} \simeq Sym_{\mathcal{O}_{X}}(\Omega^{1}_{X/k}[1]).$$
\item Si $X$ est  lisse sur $k$, alors il existe un 
isomorphisme naturel de $k$-alg\`ebres commutatives
$$\pi_{0}(\mathcal{O}(LX)^{hS^{1}}):=
\pi_{0}(\mathbb{R}\Gamma(X,S^{1}\otimes \mathcal{O}_{X})^{hS^{1}})\simeq
H^{ev}_{DR}(X/k),$$
o\`u $(-)^{hS^{1}}$ d\'esigne le foncteur des points fixes homotopiques, et 
$H^{ev}_{DR}(X/k)$ est la cohomologie de de Rham paire de $X/k$.
\end{enumerate}
\end{cor}

\textit{Preuve:} Le point $(1)$ est une cons\'equence imm\'ediate du th\'eor\`eme \ref{t1}.  
Pour le point $(2)$ il suffit de remarquer qu'il existe une \'equivalence
naturelle de faisceaux de $k$-alg\`ebres commutatives simpliciales
(o\`u le produit tensoriel d\'eriv\'e est calcul\'e dans la $sk-CAlg$)
$$S^{1}\otimes^{\mathbb{L}}\mathcal{O}_{X} \simeq \mathcal{O}_{X}\otimes^{\mathbb{L}}_{\mathcal{O}_{X}\otimes_{k}^{\mathbb{L}}
\mathcal{O}_{X}}\mathcal{O}_{X},$$
qui par normalisation donne une \'equivalence de faisceaux de cdga
(o\`u le produit tensoriel d\'eriv\'e est maintenant calcul\'e
dans $cdga$)
$$N(S^{1}\otimes^{\mathbb{L}}\mathcal{O}_{X}))\simeq 
\mathcal{O}_{X}\otimes^{\mathbb{L}}_{\mathcal{O}_{X}\otimes_{k}^{\mathbb{L}}
\mathcal{O}_{X}}\mathcal{O}_{X}.$$
D'autre part la cdga sous-jacente \`a $\mathbb{L}\epsilon(\mathcal{O}_{X})$
est naturellement \'equivalente \`a $Sym_{\mathcal{O}_{X}}(\mathbb{L}_{X/k}[1])$
comme nous l'avons d\'ej\`a fait remarqu\'e au \S 1. Le point $(3)$ se d\'eduit
formellement de $(2)$ et du fait que lorsque $X$ est lisse
$\mathbb{L}_{X/k}\simeq \Omega_{X/k}^{1}$. 

Pour le dernier point, notons $\mathcal{O}(LX):=\mathbb{R}\Gamma(X,S^{1}\otimes \mathcal{O}_{X})$, 
o\`u $\Gamma$ est le foncteur de Quillen \`a droite des sections globales, de la cat\'egorie des faisceaux 
de sur $X$ \`a valeurs dans $S^{1}-sk-CAlg$ vers la cat\'egorie $S^{1}-sk-CAlg$.
Pour $B$ une $\epsilon-cdga$ il existe un 
isomorphisme naturel d'alg\`ebres
$$Hom_{Ho(\epsilon-dg-mod)}(k,B) \simeq H^{0}(Tot(\mathcal{B}B^{-})),$$
o\`u $Tot(\mathcal{B}B^{-})$ est le complexe total n\'egatif associ\'e au complexe
mixte $B$ (voir \cite[5.1.7]{lo}, ici notre $\epsilon$ joue le r\^ole de
l'op\'erateur $B$, 
la diff\'erentielle du complexe sous-jacent \`a $B$ est $b$). En particulier, 
on a un isomorphisme naturel entre $Hom_{Ho(\epsilon-dg-mod)}(k,B)$ 
et la cohomologie du complexe de longuer $2$
$$\xymatrix{
\oplus_{n \geq 0}B_{-2n-1} \ar[r]^-{d} & \oplus_{n \geq 0}B_{-2n} \ar[r]^-{d} &
\oplus_{n \geq 0}B_{-2n+1},}$$
o\`u $d$ est la diff\'erentielle somme de $\epsilon$ et de la diff\'erentielle du complexe
sous-jacent \`a $B$. 
Ceci, appliqu\'e \`a 
$\phi(\mathcal{O}(LX))$ donne
$$\pi_{0}(\mathcal{O}(LX)^{hS^{1}}) \simeq
Hom_{Ho(S^{1}-sk-CAlg)}(k,\mathcal{O}(LX)) \simeq
Hom_{Ho(\epsilon-cdga)}(k,\phi(\mathcal{O}(LX)))\simeq
H^{0}(Tot\mathcal{B}\phi(\mathcal{O}(LX))^{-}).$$
Par $(1)$ et $(3)$ on a 
$$\phi(\mathcal{O}(LX))\simeq
\mathbb{R}\Gamma(X,\epsilon(X)))\simeq \oplus_{n}\mathbb{R}\Gamma(X,\Omega_{X}^{n}[-n]),$$
et o\`u l'action de $\epsilon$ est induite par la diff\'erentielle de de Rham sur le membre
de droite. 
Ainsi, on a donc clairement
$$H^{0}(Tot\mathcal{B}\phi(\mathcal{O}(LX))^{-})\simeq 
\oplus_{n \; pair}H^{n}_{DR}(X) = H^{ev}_{DR}(X).$$
\hfill $\Box$ \\

Il est aussi possible de g\'en\'eraliser le point $(4)$ du corollaire pr\'ec\'edent au cas non-affine
de la fa\c{c}on suivante. La cdga $\mathbb{R}\Gamma(X,S^{1}\otimes\mathcal{O}_{X})^{hS^{1}}$
est naturellement munie d'une structure de $k[u]-cdga$, et l'on peut donc
inverser $u$. Il existe alors des isomorphismes
$$\pi_{*}(\mathbb{R}\Gamma(X,S^{1}\otimes\mathcal{O}_{X})^{hS^{1}})[u^{-1}]
\simeq H^{*}_{per}(X/k),$$
o\`u $H_{per}^{i}(X/k)=H^{ev}_{DR}(X/K)$ si $i$ est pair et 
$H_{per}^{i}(X/k)=H^{odd}_{DR}(X/K)$ si $i$ est impair. Nous ne donnerons pas les d\'etails ici.

\end{document}